\documentclass[10pt]{article}

\usepackage{amsfonts}
\usepackage[ansinew]{inputenc}
\usepackage[dvips]{graphics}
\usepackage{latexsym,amssymb,amsmath}
\usepackage{graphicx}
\usepackage{exscale}
\usepackage[active]{srcltx}

\newtheorem{theorem}{Theorem}
\newtheorem{definition}{Definition}
\newtheorem{lemma}{Lemma}
\newtheorem{remark}{Remark}

\newcommand{\bq}{\begin{equation}}
\newcommand{\eq}{\end{equation}}
\newcommand{\bt}{\begin{tabular}}
\newcommand{\et}{\end{tabular}}
\newcommand{\ba}{\begin{array}}
\newcommand{\ea}{\end{array}}
\newcommand{\be}{\begin{equation}}
\newcommand{\ee}{\end{equation}}
\newcommand{\bes}{\begin{equation*}}
\newcommand{\ees}{\end{equation*}}

\begin{document}

\date{}
\title{First order non-homogeneous $q$-difference equation for Stieltjes
function characterizing $q$-orthogonal polynomials}
\author{J. Arves\'u, and A. Soria-Lorente \\
Department of Mathematics,\\
Universidad Carlos III de Madrid,\\
Avda. de la Universidad, 30, 28911, Legan\'es, Madrid, Spain}
\date{\today }
\maketitle

\begin{abstract}
In this paper we give a characterization of some classical $q$-orthogonal polynomials in terms of a difference
property of the associated \emph{Stieltjes function}, i.e this function solves a
first order non-homogeneous $q$-difference equation. The solutions of the
aforementioned $q$-difference equation (given in terms of hypergeometric
series) for some canonical cases, namely, \emph{$q$-Charlier}, 
\emph{$q$-Kravchuk}, \emph{$q$-Meixner} and \emph{$q$-Hahn} are worked out.
\end{abstract}




\section{Introduction}

\noindent The goal of this paper is to extend to the $q$-discrete orthogonal
polynomials on the non-uniform lattice $x(s)=\left( q^{s}-1\right) /\left(
q-1\right) $, the characterization for classical orthogonal polynomials in
terms of the standard Stieltjes function (see \cite{stieltjes} for original
and historical issues on the Stieltjes function; for modern and advanced
notions see \cite{Chihara}, \cite{ma-la}, \cite{simon} and the references
therein).

The properties of $q$-\emph{Charlier}, $q$-\emph{Kravchuk}, $q$-\emph{Meixner} and $q$-\emph{Hahn} orthogonal families are well known (see \cite
{Arv1,Niki} and \cite{me-al-ma}): they satisfy a hy\-per\-geo\-me\-tric-type
difference equation (\ref{eqdif}); their finite differences constitute an
orthogonal polynomial family (\ref{ortho-diffs}); they can be expressed by a
Rodrigues-type formula (\ref{rodrigues-t}); their associated orthogonalizing
weights satisfy a Pearson-type difference equation (\ref{pearson})-(\ref
{pearson2}); they verify a three-term recurrence relation (\ref{ttrrc}).
These properties, among others, characterize the classical orthogonal
polynomials (see e.g. \cite{Al,ma-bra-pet,ma-pet,maroni2}). Most of these characterization properties have been already extended to the $q$-polynomials. However, special attention must be given to \cite{me-al-ma}, particularly for the unified study of some classical $q$-polynomials carried out. An algebraic approach developed by Maroni \cite{maroni} is used as a crucial ingredient. Indeed, the starting point is a distributional difference equation fulfilled by the moment functional with respect to which the aforementioned $q$-polynomials are orthogonal. The motivation for studying the problem considered in the present paper comes precisely from an unproven assertion given in Proposition 2.27 (b) of the aforementioned paper \cite{me-al-ma}. In this paper we will prove the above mentioned assertion for classical $q$-polynomials on the lattice $(q^s-1)/(q-1)$, i.e., a $q$-moment functional $\mathcal{U}$ is classic if and only if the Stieltjes function associated with it, given in terms of the $q$-falling factorial, solves a non-homogeneous version of the Pearson-type difference equation associated with $\mathcal{U}$ (see Theorem \ref{main-th} below).

The advantages of the proposed approach over existing works consist of
working directly with linear functionals rather than power series expansion
for the Stieltjes function, i.e., of looking for the intrinsic properties
that may satisfy the linear functional under consideration acting, either on
the vector space of formal power series or on the vector space of
polynomials. Our characterization provides a general framework which can be
extended -like in the continuous case- to the study of $q$-semiclassical
orthogonal polynomials \cite{me-al-ma}. In addition, for illustrating our
main result, the Stieltjes function given in terms of the $q$-falling
factorial for the four canonical cases, namely, \emph{q-Charlier}, \emph{\
q-Kravchuk}, \emph{q-Meixner}, and \emph{q-Hahn} is explicitly calculated
and expressed in terms of hypergeometric functions. As far as we know, there
were not explicit hypergeometric expressions for these Stieltjes functions
in the literature.

The structure of the paper is as follows. In Section \ref{pre_notions}, we
give a quick overview of the relationship between orthogonal polynomials and
Stieltjes function. In Section \ref{def-notations}, we compress some
necessary definitions and tools. Lastly, in Section \ref{q-Stieltjes} the
Stieltjes function for $q$-orthogonal polynomials in terms of the $q$
-falling factorial basis is introduced and the main theorem for the proposed
characterization is proved. This result constitutes a $q$-analogue of the
characterization given in \cite{ma-la} for classical discrete orthogonal
polynomials (see also \cite{ga-ma-la}). Finally, in Section \ref{examples}
the explicit expressions for the corresponding solutions of the difference
equation that characterizes the classical $q$-orthogonal polynomials on the
non-uniform lattice $x(s)=(q^{s}-1)/(q-1)$, i.e. the Stieltjes functions,
are studied.

\section{Stieltjes function and orthogonal polynomials\label{pre_notions}}

We start by reviewing some results needed for the sequel; mainly extracted from \cite{Chihara,Niki,Nikishin}. The linear space
of polynomials with coefficients in $\mathbb{C}$ (the set of complex
numbers) is denoted by $\mathbb{P}$, and its dual space by $\mathbb{P}^{*}$,
whose elements are linear functionals (moment functionals).

Let the action of $\mathcal{U}\in \mathbb{P}^{\ast }$ on $p\in \mathbb{P}$
be denoted by the duality bracket $\left\langle \mathcal{U},p\right\rangle$. The moments of $\mathcal{U}$ are given by $\left\langle \mathcal{U},x^{k}\right\rangle =u_{k}$,
$k=0,1,\ldots$, and by linearity yields the following relation $\left\langle \mathcal{U},P_{n}(x)\right\rangle
=\sum_{k=0}^{n}\alpha _{k}u_{k}$, where $P_{n}(z)=\sum_{k=0}^{n}\alpha _{k}z^{k}$.

Formally, if $\mathcal{U}$ is a moment functional determined by the moment
sequence $\{u_{k}\}_{k\geq 0}$, the Stieltjes function associated with $
\mathcal{U}$ can be defined as the generating function of these moments by
means of the series expansion 
\begin{equation}
S(z)=\sum_{k\geq 0}\frac{u_{k}}{z^{k+1}}.  \label{ref1}
\end{equation}
Notice that $S(z)$ corresponds to its sum when the Laurent expansion (\ref
{ref1}) converges for some $z$; otherwise $S(z)$ represents its analytic
continuation -provided that it exists-.

The above function (\ref{ref1}) can be approximated by rational functions
with prescribed order near infinity \cite{Nikishin}. This approximation
problem -often known as $(n-1,n)$ Pad\'e approximation near infinity- consists
in finding the polynomials $P_{n}(z)=\sum_{k=0}^{n}\alpha _{k}z^{k}$ and $
Q_{n-1}(z)=\sum_{k=0}^{n-1}\beta _{k}z^{k}$, such that 
\begin{equation}
P_{n}(z)S(z)-Q_{n-1}(z)=\mathcal{O}(z^{-n-1}),  \label{pade}
\end{equation}
where $\deg P_{n}\leq n$ and $\deg Q_{n-1}\leq n-1$.

By getting $Q_{n-1}(z)$ equal to the polynomial part of $P_{n}(z)S(z)$ and
using the prescribed order at infinity in the interpolation condition (\ref{pade}), the
following linear system of $n$ homogeneous equations \cite{Nikishin} 
\begin{equation}
\displaystyle\sum_{i=0}^{n}u_{i+j}\alpha _{i}=0,\quad j=0,\dots ,n-1,
\label{system-eq-pade}
\end{equation}
for determining the $n+1$ unknown coefficients of $P_{n}(z)$, holds. Hence, $
P_{n}(z)$ is determined up to a multiplicative factor. Observe that the linear system (\ref{system-eq-pade}) is equivalent to the following orthogonality conditions: 
\begin{equation}
\left\langle \mathcal{U},P_{n}(x)x^{k}\right\rangle =0,\quad k=0,\dots ,n-1.
\label{ortho}
\end{equation}
Consequently, $P_{n}(z)$ is orthogonal with respect to any polynomial of degree less than $n$ for the linear functional $\mathcal{U}$. On the other hand, we have selected $Q_{n-1}(z)$ equals to the polynomial part of $P_{n}(z)S(z)$, hence one has
\begin{equation}
\displaystyle Q_{n-1}(z)=\left\langle \mathcal{U},\frac{P_{n}(z)-P_{n}(x)}{
z-x}\right\rangle,  \label{Q-pade}
\end{equation}
where it is assumed that $\mathcal{U}$ acts on the variable $x$. Notice that (\ref{Q-pade}) leads to the formal relation
\begin{equation}
P_{n}(z)\left\langle \mathcal{U},\frac{1}{(z-x)}\right\rangle-Q_{n-1}(z)=\left\langle \mathcal{U},\frac{P_{n}(x)}{z-x}\right\rangle,
\label{z-x}
\end{equation}
which indeed implies -at least formal- the linearity extension of the linear functional $\mathcal{U}$ to the vector space of formal power series. Moreover, regardless of the convergence of (\ref{ref1}), it is straightforward to see that this series provides a means of approximating $\left\langle \mathcal{U},(z-x)^{-1}\right\rangle$. Indeed, for every nonnegative integer $n$,
$$
\left\langle\mathcal{U},\dfrac{1}{z-x}\right\rangle =
\left\langle\mathcal{U},\sum_{k=0}^{n-1}\dfrac{x^{k}}{z^{k+1}}\right\rangle
+\left\langle\mathcal{U},\left(\dfrac{x}{z}\right)^{n}\dfrac{1}{z-x}\right\rangle.$$
Thus, for $z\to\infty$, in every sector
$\epsilon<\mbox{arg}z<\pi-\epsilon$ $(0<\epsilon<\pi/2)$ one gets $\vert z-x\vert\geq\vert z\vert\sin\epsilon$. Hence,
$$
\left\langle\mathcal{U},\dfrac{1}{z-x}\right\rangle -\sum_{k=0}^{n-1}\frac{u_{k}}{z^{k+1}}=
\mathcal{O}_{\epsilon}\left(z^{-n-1}\right),
$$
which holds for all $z$ in the above sector. Analogously, for the remainder approximation term we get the prescribed order $z^{-n-1}$ near $\infty$, in the above sector, i.e
\begin{equation*}
\left\langle \mathcal{U},\frac{P_{n}(x)}{z-x}\right\rangle =\mathcal{O}_{\epsilon}
(z^{-n-1}).
\end{equation*}
Moreover, we recall that in the framework of the algebraic approach developed by Maroni \cite{maroni,maroni2,maroni3} it is established a topological isomorphism between the space $\mathbb{P}^{*}$ and the space of formal power series endowed with appropriate topologies. Indeed, it is always possible to associate $(P_n)_{n\geq0}$ of degree $n$ with a unique sequence $(\mathcal{U}_{n})_{n\geq0}$, $\mathcal{U}_{n}\in\mathbb{P}^{*}$, called the dual sequence of $(P_n)_{n\geq0}$, such that $\left\langle\mathcal{U}_{n},P_m\right\rangle=0$, for $n\not=m$, and equals $1$ for $n=m$. Indeed, any functional $\mathcal{U}\in\mathbb{P}^{*}$ can be represented as follows
$$
\mathcal{U}=\sum_{n\geq0}\left\langle\mathcal{U},P_n\right\rangle\mathcal{U}_{n},
$$
which makes use of the linearly extension of $\mathcal{U}$ to the vector space of formal power series. Therefore, the Stieltjes function associated with $\mathcal{U}$ in the Pad\'e approximation problem (\ref{pade}) has
the formal representation 
\begin{equation}
S(z)=\left\langle \mathcal{U},(z-x)^{-1}\right\rangle.  \label{ref2-formal}
\end{equation}

From now on we will consider a specific class of moment functionals. Let $\mu$ be a Borel measure on the real line 
$\mathbb{R}$ (with infinitely many points of increase), supported on 
$\Omega \subset \mathbb{R}$. We define the linear functional on which we will focus our attention
\begin{equation}
\left\langle \mathcal{U},p\right\rangle =\int_{\Omega
}p(x)\,d\mu (x),\quad p\in\mathbb{P}.\label{orthoU}
\end{equation}
In particular, we will consider
\begin{equation*}
\mu =\sum_{k=0}^{N}\rho (x)\delta _{x_{k}},\quad \rho (x_{k})>0,\ x_{k}\in 
\mathbb{R},\quad N\in \mathbb{N}\cup \{+\infty \}.
\end{equation*}
This measure $\mu$ is discrete (with finite moments $u_{k}=\int_{\Omega }x^{k}d\mu (x)$, $k\geq0$) formed by a linear combination of Dirac measures at the points $x_{0},\ldots
,x_{N}$. By $\mathbb{N}$ we denotes the set of all nonnegative integers. Therefore, the orthogonality condition (\ref{ortho}) for the $q$-polynomials on the lattice $\{x(s)\mapsto 
\mathbb{R}^{+}:\,s=0,\dots,N\}$ with respect to (\ref{orthoU}) is defined  as follows --see \cite[eq. (4.7), p.256]{slov} as well as \cite{nubook} 
\begin{equation}
\displaystyle\left<\mathcal{U},P_{n}x^k\right>=\sum
\limits_{s=0}^{N}P_n(x(s))x^{k}(s)\,\rho(s)\bigtriangleup x(s- 1/2 )=0, \quad k=0,\dots,n-1,  \label{d-ortho-c}
\end{equation}
where, $\bigtriangleup x(s)=x(s+1)-x(s)$ denotes the forward difference
operator. In general, the polynomial $P_n(x(s))$ is called discrete orthogonal
polynomial. In what follows we will denote any polynomial $P_n(x(s))$ simply as $
P_n(s)$.

The term \emph{`classical discrete orthogonal polynomial'} \cite{Niki} (also
known as classical orthogonal polynomials of a discrete variable) obeys the fact that $x(s)=c_{1}q^{s}+c_{2}q^{-s}+c_{3}$, $(q\in\mathbb{R}^{+}\backslash{\{1\}})$
or $x(s)=c_4s^{2}+c_5s+c_6$, where the constants $c_{i}\in\mathbb{R}$ $
(i=1,\dots,6)$ are independents of variable $s$, and the orthogonalizing
weights $\rho(s)$ are solutions of the Pearson-type difference equation (see
formulas (\ref{pearson})-(\ref{pearson2}) below and \cite[pp. 70-72]{Niki}).

A remarkable feature of the classical discrete orthogonal polynomials is that they satisfy the hypergeometric-type difference equation 
\begin{equation}
\begin{array}{c}
\displaystyle\sigma (s)\frac{\bigtriangleup }{\bigtriangleup x(s- 
1/2)}\frac{\bigtriangledown y(s)}{
\bigtriangledown x(s)}+\tau (s)\frac{\bigtriangleup y(s)}{\bigtriangleup
x(s) }+\lambda y(s)=0, \\ 
\sigma (s)=\tilde{\sigma}(x(s))-(\tilde{\tau}
(x(s))\bigtriangleup x(s-1/2))/2,\quad \tau (s)= 
\tilde{\tau}(x(s)),
\end{array}
\label{eqdif}
\end{equation}
where $\bigtriangledown y(s)=\bigtriangleup y(s-1)$ is the backward difference
operator, and 
\begin{align}
\widetilde{\sigma }(s)& =a_{2}x^{2}\left( s\right) +a_{1}x\left( s\right)
+a_{0}  \label{sigma-coeff} \\
\widetilde{\tau }(s)& =b_{1}x(s)\ +b_{0}.  \label{tau-coeff}
\end{align}
Equation (\ref{eqdif}) is a discretization of the very well known
hypergeometric differential equation (see \cite{Niki}) 
\begin{equation*}
\displaystyle\tilde{\sigma}(x)y^{\prime \prime }(x)+\tilde{\tau}(x)+\tilde{
\lambda}y(x)=0,\quad \mbox{with}\,\,\deg \tilde{\sigma}\leq 2,\quad \deg 
\tilde{\tau}=1,\quad \tilde{\lambda}\in \mathbb{R}.
\end{equation*}
Notice that the solution $\rho (s)$ of the Person-type equation (written in
two equivalent forms) 
\begin{align}
\displaystyle\Delta \left[ \sigma (s)\rho (s)\right] & =\tau (s)\rho
(s),\quad \Delta \overset{\mbox{\tiny def}}{=}\frac{\bigtriangleup }{
\bigtriangleup x(s-1/2)},  \label{pearson} \\
\nabla \lbrack (\sigma (s)+\tau (s)\bigtriangledown x(s+1/2))\rho (s)]& =\tau (s)\rho (s),\quad \nabla \overset{ 
\mbox{\tiny def}}{=}\frac{\bigtriangledown }{\bigtriangledown x(s+ 1/2)},  \label{pearson2}
\end{align}
allows to write (\ref{eqdif}) in self-adjoint form. This symmetrization
factor of (\ref{eqdif}) with the condition $\displaystyle\left. \sigma
(s)\rho (s)x^{k}(s-1/2)\right\vert _{s=0,N+1}=0$, $k=0,1,\dots $, is called the orthogonalizing weight of (\ref{d-ortho-c}).
This boundary condition is compatible with the existence of all moments of 
(\ref{d-ortho-c}).

It is well-known that many properties of orthogonal polynomials can be
deduced from (\ref{d-ortho-c}); in particular, equation (\ref
{eqdif}). In turn from this equation follow important properties like the
Rodrigues-type formula (see \cite{Niki,slov} for general theory and \cite{Arv1} for main data of $q$-polynomials on non-uniform lattice $x(s)=(q^s-1)/(q-1)$). 
\begin{equation}
P_{n}(s)=\frac{B_{n}}{\rho (s)}\nabla ^{(n)}[\rho _{n}(s)],\quad \nabla
^{(n)}=\frac{\bigtriangledown }{\bigtriangledown x_{1}(s)}\frac{
\bigtriangledown }{\bigtriangledown x_{2}(s)}\cdots \frac{\bigtriangledown }{
\bigtriangledown x_{n}(s)},  \label{rodrigues-t}
\end{equation}
where $x_{k}(s)=x(s+\frac{k}{2})$, and $\rho_n(s) = \rho(s+n) \prod_{m=1}^n
\sigma(s+m)$; the orthogonality of the differences 
\begin{equation}
\sum_{s=0}^{N-1}\frac{\bigtriangleup P_{n}(s)}{\bigtriangleup x(s)}\frac{
\bigtriangleup P_{m}(s)}{\bigtriangleup x(s)}\rho _{1}(s)\bigtriangleup
x(s)=\delta _{n,m}\left\Vert \frac{\bigtriangleup P_{n}(s)}{\bigtriangleup
x(s)}\right\Vert ^{2},  \label{ortho-diffs}
\end{equation}
where $\left\Vert\cdot\right\Vert$ denotes the norm in the Hilbert space $L^2_\mu$; the three-term recurrence relation 
\begin{equation}
x(s)P_{n}(s)=\alpha _{n}P_{n+1}(s)+\beta _{n}P_{n}(s)+\gamma
_{n}P_{n-1}(s),\quad n=0,1,\dots ,N,  \label{ttrrc}
\end{equation}
where $P_{-1}=0$, $P_{0}$ is a real constant, and the coefficients $\alpha
_{n}$, $\beta _{n}$, and $\gamma _{n}$ are given explicitly in terms of 
(\ref{sigma-coeff}) and (\ref{tau-coeff}) (see \cite{Arv1}, expressions in (42)).

A second linearly independent solution of equation (\ref{eqdif}) --for a
specific choice $\lambda=\lambda_n$, see \cite{Niki}-- is the so-called
function of the second kind on non-uniform lattice (see \cite{slov}, formula
(5.1)). This function is also a second linearly independent solution of
equation (\ref{ttrrc}). For a relationship between the function of the
second kind and the Stieltjes function we refer to \cite{Niki,slov}.

\section{Basic definitions and notations\label{def-notations}}

Here, we will use the $q$-analogue of the Pochhammer symbol \cite
{Gasper,koekoek,Niki} 
\begin{equation}
(a;q)_{k}=\prod_{0\leq j\leq k-1}(1-aq^{j}),\quad\mbox{for}\quad k>0,\quad 
\text{and}\quad (a;q)_{0}=1.  \label{q-pochhammer}
\end{equation}
By Pochhammer symbol $(s)_{k}$ and falling factorial $[s] _{k}$ we mean 
\begin{align*}
(s)_{k}&=s(s+1)\cdots (s+k-1),\quad (s)_{0}=1,\quad k\geq 1, \\
\left[ s\right] _{k}&=(-1)^{k}(-s)_{k},
\end{align*}
respectively. Notice that the subscript in the above expressions is a nonnegative integer, whereas in the sequel symbol $\left[s\right] _{q}$ denotes the $q$-number defined as follows
\begin{equation*}
\left[ s\right] _{q}=\frac{q^{s/2}-q^{-s/2}}{q^{1/2}-q^{-1/2}},\quad s\in \mathbb{C}.
\end{equation*}
For the $q$-factorial we use the following definition
\begin{equation*}
\left[ n\right] _{q}!=\left[ n\right] _{q}\left[ n-1\right] _{q}\cdots \left[
2\right] _{q}\left[ 1\right] _{q},\quad n\in \mathbb{N},\quad 
\left[ 0\right] _{q}!=1.
\end{equation*}
Observe that in \cite[p. 7]{Gasper} the $q$-number is defined in a different way. 

The $q$-hypergeometric series $_{r}\varphi _{s}$ is defined as 
\begin{equation*}
_{r}\varphi _{s}\left( 
\begin{array}{c|c}
a_{1},\ldots ,a_{r} &  \\ 
& q;z \\ 
b_{1},\ldots ,b_{s} & 
\end{array}
\right) =\sum_{k\geq 0}\frac{\left( a_{1};q\right) _{k}\cdots \left(
a_{r};q\right) _{k}}{\left( b_{1};q\right) _{k}\cdots \left( b_{s};q\right)
_{k}}\left[ \left( -1\right) ^{k}q^{\binom{k}{2}}\right] ^{1+s-r}\frac{z^{k} 
}{\left( q;q\right) _{k}},
\end{equation*}
where $\binom{m}{n}$ denotes the binomial coefficient.

In what follows, the $q$-falling factorial defined as 
\begin{equation}
\left[ s\right] _{q}^{(k)}=\prod_{0\leq j\leq k-1}x(s-j),\quad \mbox{for}
\quad k>0,\quad \text{and}\quad \left[ s\right] _{q}^{(0)}=1,
\label{q-falling-basis}
\end{equation}
where $x(s)=\frac{q^{s}-1}{q-1}$, will play the same role that the power $
x^{k}$ plays in a formal series expansion. Indeed, the $q$-falling factorial
basis $\{1,\left[ s\right] _{q}^{(1)},\ldots ,\left[ s\right]
_{q}^{(n)},\dots \}$ has been suggested to be more natural than the basis $
\{1,x(s),\dots ,x^{n}(s),\dots \}$, mainly when multiple orthogonal
polynomials are considered \cite{arvesu-q-multi-Hahn}. Observe that the $q$
-falling factorial $\left[ s\right] _{q}^{(k)}$ is a polynomial of degree $k$
in $x(s)$. Accordingly, the action on it of the forward and backward
difference operators defined in (\ref{pearson})-(\ref{pearson2}) gives 
\begin{equation*}
\Delta \left[ s\right] _{q}^{(k)}=q^{3/2-k}\left[ k\right] _{q}^{(1)}\left[
s \right] _{q}^{(k-1)},\,\,\mbox{and}\,\,\nabla \left[ s\right]
_{q}^{(k)}=q^{1/2-k}\left[ k\right] _{q}^{(1)}\left[ s-1\right] _{q}^{(k-1)},
\end{equation*}
respectively.

In addition, the $q$-falling factorial can be rewritten in terms of the $q$
-analogue of the Pochhammer symbol as follows 
\begin{equation}
\left[ s\right] _{q}^{(k)}=\frac{(q^{-s};q)_{k}}{(q-1)^{k}}q^{k\left( s-
\frac{k-1}{2}\right) },\quad k\geq 1,\quad \mbox{and}\quad \left[ s\right]
_{q}^{(0)}=1.  \label{ref-q-falling}
\end{equation}
Clearly, the Pochhammer symbol and the falling factorial can be recovered as
a limiting case 
\begin{equation*}
\lim_{q\rightarrow 1}\left[ s\right] _{q}^{(k)}=(-1)^{k}(-s)_{k}=\left[ s
\right] _{k}.
\end{equation*}
Let $\mathbb{S}$ be the spanning set of $\{1,\left[ s\right]
_{q}^{(1)},\ldots ,\left[ s\right] _{q}^{(n)},\dots \}$ over $\mathbb{R}$.
Observe that $\mathbb{S}$ coincides with the linear space of polynomials of
discrete variable $x(s)$ with real coefficients. $\mathbb{S}_{n}$ denotes
the corresponding linear subspace of dimension $n$. Let us define the
algebraic dual space of $\mathbb{S}$, denoted by $\mathbb{S}^{\ast }$, as 
\begin{equation*}
\mathbb{S}^{\ast }=\{\mathcal{U}:\mathbb{S}\rightarrow \mathbb{R},\hspace{
0.1in}\mbox{such that $\mathcal{U}$ is a linear
functional}\}.
\end{equation*}
The real number $u_{n}^{q}=\left\langle \mathcal{U},[s]_{q}^{(n)}\right
\rangle $, $n\geq 0$, \label{q-moment-U} is said to be the $q$-moment of
order $n$ of $\mathcal{U}$ and the sequence $\{u_{n}^{q}\}_{n\geq 0}$ is
called the $q$-moment sequence associated with $\mathcal{U}$.

\begin{definition}
\label{first-def} Let $\{u_n^q \}_{n\geq 0}$ be a sequence in $\mathbb{R}$. $
\mathcal{U}$ is said to be the $q$-moment functional associated with the
sequence $\{u_n^q \}_{n\geq 0}$, if $\langle\mathcal{U},[s]_q^{(n)}
\rangle=u_n^q$ and it can be extended to $\mathbb{S}$ by linearity, i.e., if 
\begin{equation*}
\displaystyle r(s)=\sum_{k=0}^{n}a_k[s]_q^{(k)},\,\,\mbox{then}
\,\,\left\langle\mathcal{U},r\right\rangle=\sum_{k=0}^{n}a_k\mathit{u}_k^q.
\end{equation*}
\end{definition}

Accordingly, from Definition \ref{first-def} and equation (\ref{d-ortho-c})
one has 
\begin{equation}
\langle \mathcal{U},r\rangle =\sum_{s=0}^{N}r(s)\rho (s)\bigtriangleup x(s- 
\mbox{\scriptsize $\frac{1}{2}$}),  \label{def-q-functional}
\end{equation}
where, the $q$-moment of order $k$ is given by 
\begin{equation}
u_{k}^{q}=\left\langle \mathcal{U},\left[ s\right] _{q}^{(k)}\right\rangle
=\sum_{s=0}^{N}\left[ s\right] _{q}^{(k)}\rho (s)\bigtriangleup x(s- 
\mbox{\scriptsize $\frac{1}{2}$}),\quad \quad \quad k=0,1,\ldots
\label{q-moment-u}
\end{equation}
Below we sketch some definitions and results that we will use in the next
sections (for comprehensive definitions and theorems on $q$-moment
functionals we refer to \cite{me-al-ma}).

Any moment functional $\mathcal{U}$, which maps $\mathbb{S}$ into $\mathbb{S}
$ has a transpose $\mathcal{V}^{\mathrm{T}}$ : $\mathbb{S}^{*}\mapsto\mathbb{
\ \ S}^{*}$ defined by $\langle\mathcal{V}^{\mathrm{T}}\mathcal{U}
,p\rangle=\langle\mathcal{U},\mathcal{V}p\rangle$, $\forall p\in\mathbb{S}$, 
$\forall\mathcal{U}\in\mathbb{S}^{*}$; then $\Delta \mathcal{U}$, $\nabla 
\mathcal{U}$, and $p\mathcal{U}$ are defined by duality according to the
following definitions.

\begin{definition}
\label{def-delta-nabla-U} Let $\mathcal{U}$ be a $q$-moment functional. The
forward and backward differences of the moment functional $\mathcal{U}$,
denoted by $\Delta \mathcal{U}$, and $\nabla\mathcal{U}$, respectively, are
defined as follows 
\begin{align*}
\langle \Delta \mathcal{U},p(s)\rangle &=-\langle \mathcal{U},\nabla
p(s)\rangle,\quad\mbox{where}\quad p\in\mathbb{S}, \\
\langle \nabla \mathcal{U},p(s)\rangle &=-\langle \mathcal{U},\Delta
p(s)\rangle.
\end{align*}
\end{definition}

Here, the forward difference operator $\Delta$ on moment functionals is
minus the transpose of the backward difference operator on polynomials.

\begin{definition}
\label{def-prod-pU} Let $p \in \mathbb{S}$. The linear functional $\mathcal{
\ \ V }=p\mathcal{U}$ is said to be the left-multiplication of $\mathcal{U}$
by a polynomial $p$ if 
\begin{equation*}
\left\langle \mathcal{V},r\right\rangle=\left\langle\mathcal{U}
,pr\right\rangle, \quad r\in \mathbb{S}.
\end{equation*}
\end{definition}

\begin{definition}
If $\left( \sigma ,\tau \right)$ are polynomials of minimum degree such that
\begin{equation}
\Delta (\sigma \mathcal{U})=\tau \mathcal{U},
\label{classic-func}
\end{equation}
define the class of $\mathcal{U}$ as the nonnegative integer number $s$ such that
$$
s = \max\{\deg \sigma-2, \deg \tau-1\}.
$$
The polynomial sequence $\left\{ P_{n}(s)\right\} _{n\geq 0}$ orthogonal with
respect to the $q$-moment functional $\mathcal{U}$ of class zero is said to be classical orthogonal polynomial sequence.
\end{definition}

Equation (\ref{classic-func}) is a $q$-analogue of the functional Person
equation \cite{garri-arvesu-paco,ma-la}. In general, the polynomial
sequences orthogonal with respect to such moment functionals $\mathcal{U}$
are said to be semiclassical orthogonal polynomial sequences.

\section{Stieltjes function for $q$-orthogonal polynomials in terms of the $
q $-falling factorials \label{q-Stieltjes}}

Assume that $\mathcal{U}$ is the moment functional determined by the $q$
-moment sequence $\{u^{q}_k\}_{k\geq0}$. Below, in the same fashion that the
continuous case (\ref{ref1}) we will define a $q$-analogue of the Stieltjes
function associated with $\mathcal{U}$ (in short $q$-Stieltjes function) as
the generating function of these $q$-moments by means of a series expansion,
but using the $q$-falling factorial basis (\ref{q-falling-basis}) instead.
More precisely,

\begin{definition}
Let $\mathcal{U}$ be a $q$-moment functional defined on the linear space $
\mathbb{S}$. The $q$-Stieltjes function associated with $\mathcal{U}$ is
defined as the formal series expansion 
\begin{equation}
S_{q}(z)=\sum_{k\geq 0}\frac{u_{k}^{q}}{q^{k}\left[ z\right] _{q}^{(k+1)}}.
\label{q-stieltjes}
\end{equation}
\end{definition}

The following result deals with the $q$-analogue of the formal representation
(\ref{ref2-formal}). Indeed, we will use a $q$-analogue of the
Chu-Vandermonde convolution (see \cite{koekoek} for details) 
\begin{equation*}
_{2}\varphi _{1}\left( 
\begin{array}{c|c}
q^{-n},b &  \\ 
& q;q^{n}c/b \\ 
c & 
\end{array}
\right) =\frac{(cb^{-1};q)_{n}}{(c;q)_{n}}.
\end{equation*}
In particular, one can select $s\in \{0,1,\dots ,N\}$ and $z\in \mathbb{C}$ such that
\begin{equation}
_{2}\varphi _{1}\left( 
\begin{array}{c|c}
q^{-s},q &  \\ 
& q;q^{s-z} \\ 
q^{1-z} & 
\end{array}
\right) =\frac{(q^{-z};q)_{s}}{(q^{1-z};q)_{s}}.  \label{q-Chu-Vander}
\end{equation}
We will use this expression in the proof of the following lemma as well as the linearity extension of $\mathcal{U}$ to the linear space of formal power series in accordance with the discussion given in section \ref{pre_notions} with the usual adaptation required for discrete measures. Indeed, using the relation (see p. 262 \cite{slov})
$$
\dfrac{1}{x(s)-x(z)}=-\dfrac{1}{x(z)}\sum_{k=0}^{n-1}\left(\dfrac{x(s)}{x(z)}\right)^{k}+\dfrac{x^{n}(s)}{[x(s)-x(z)]x^{n}(z)},
$$
and proceeding like in (\ref{z-x}) one gets an asymptotic expansion of (\ref{q-stieltjes-Suslov}) as $x(z)\to\infty$. However, we refer to \cite{maroni,maroni2,maroni3} for a more comprehensive analysis.
\begin{lemma}
\label{lemma-first}The $q$-Stieltjes function (\ref{q-stieltjes}) admits the
following formal representation 
\begin{equation}
S_{q}(z)=\left\langle \mathcal{U},\left[ x(z)-x(s)\right] ^{-1}\right\rangle
,  \label{formal-representation}
\end{equation}
where it is assumed that $\mathcal{U}$ acts on the discrete variable $x(s)$.
\end{lemma}

\emph{Proof}. Relations (\ref{ref-q-falling}) and (\ref{q-moment-u}) enable
us to rewrite (\ref{q-stieltjes}) as follows

\begin{equation}
\displaystyle S_{q}(z)=\sum_{k\geq 0}\frac{\left\langle \mathcal{U},\left[ s 
\right] _{q}^{(k)}\right\rangle }{q^{k}\left[ z\right] _{q}^{(k+1)}}= 
\displaystyle\sum_{k\geq 0}\frac{(q-1)\left\langle \mathcal{U}
,(q^{-s};q)_{k}q^{ks-\binom{k}{2}}\right\rangle }{
q^{k}(q^{-z};q)_{k+1}q^{(k+1)z-\binom{k+1}{2}}}.  \label{Sq}
\end{equation}
Particularizing the well known relation $\left( a;q\right) _{k+n}=\left(
a;q\right) _{n}\left( aq^{n};q\right) _{k}$ -see \cite{Gasper}- for the
choice $a=q^{-z}$ and $n=1$, i.e. 
\begin{equation*}
\left( q^{-z};q\right) _{k+1}=\left( q^{-z};q\right) _{1}\left(
q^{1-z};q\right) _{k},
\end{equation*}
it follows that equation (\ref{Sq}) can be rewritten as 
\begin{align}
S_{q}(z)& =\frac{q-1}{q^{z}-1}\sum_{k\geq 0}\frac{\left\langle \mathcal{U}
,q^{k(s-z)}(q^{-s};q)_{k}\right\rangle }{(q^{1-z};q)_{k}}  \notag \\
& =\frac{1}{x\left( z\right) }\left\langle \mathcal{U},\sum_{k\geq 0}\frac{
(q^{-s};q)_{k}(q;q)_{k}}{(q^{1-z};q)_{k}(q;q)_{k}}\Big(\frac{q^{s}q^{1-z}}{q}
\Big)^{k}\right\rangle .  \label{pre-Chu}
\end{align}
If one uses (\ref{q-Chu-Vander}) equation (\ref{pre-Chu}) transforms into
the expression 
\begin{equation}
S_{q}(z)=\left\langle \mathcal{U},\frac{1}{x\left( z\right) }\frac{
(q^{-z};q)_{s}}{(q^{1-z};q)_{s}}\right\rangle .  \label{sq-final}
\end{equation}
Taking into account relation (\ref{q-pochhammer}), a straightforward
calculation leads to the following equalities 
\begin{equation*}
\displaystyle\frac{1}{x\left( z\right) }\frac{(q^{-z};q)_{s}}{
(q^{1-z};q)_{s} }=\frac{(1-q)(1-q^{-z})}{(1-q^{z})(1-q^{-z+s})}=\frac{1}{
x(z)-x(s)},
\end{equation*}
which makes expression (\ref{sq-final}) coincident with (\ref
{formal-representation}); therefore the statement holds.

\begin{remark}
From (\ref{def-q-functional}) and (\ref{formal-representation}) one can
write (\ref{q-stieltjes}) up to a constant factor as 
\begin{equation}
S_{q}(z)=\sum_{s=0}^{N}\frac{\rho (s)\bigtriangleup x(s-1/2)}{x(z)-x(s)},  \label{q-stieltjes-Suslov}
\end{equation}
which reveals its singularities.
\end{remark}

The Stieltjes function (\ref{q-stieltjes-Suslov}) is related to the function
of the second kind on non-uniform lattice (see \cite{slov}, formula (5.1))
for the specific choice $n=0$. This function of the second kind is the
second linearly independent solution of equations (\ref{eqdif}) and (\ref
{ttrrc}), respectively.

In the next results the subscripts are used to indicate on which variable
the difference operators are acting.

\begin{lemma}
\label{diff-1-Stieltjes} The following relation 
\begin{equation}
\left\langle \nabla _{s}\mathcal{U},\frac{1}{x(z)-x(s)}\right\rangle
=\left\langle \mathcal{U}\,,\nabla _{z}\frac{1}{x(z)-x(s)}\right\rangle ,
\label{ext-def-2}
\end{equation}
holds.
\end{lemma}

\emph{Proof}. By using the operator $\nabla $ defined in (\ref{pearson2})
one can compute the backward difference of the Stieltjes function given in
the $q$-falling factorial basis. Indeed, by Lemma \ref{lemma-first} one gets 
\begin{align*}
\nabla S_{q}(z)& =\frac{1}{\bigtriangledown x(z+ 
\mbox{\scriptsize
$\frac{1}{2}$})}\left[ S_{q}(z)-S_{q}(z-1)\right] \\
& =\frac{1}{\bigtriangledown x(z+\mbox{\scriptsize $\frac{1}{2}$})}\left[
\left\langle \mathcal{U},\frac{1}{x(z)-x(s)}\right\rangle -\left\langle 
\mathcal{U},\frac{1}{x(z-1)-x(s)}\right\rangle \right] ,
\end{align*}
or equivalently, 
\begin{equation*}
\nabla S_{q}(z)=\left\langle \mathcal{U},\nabla _{z}\frac{1}{x(z)-x(s)}
\right\rangle .
\end{equation*}
A straightforward computation leads to the following relation 
\begin{equation*}
\nabla _{z}\frac{1}{x(s)-x(z)}=\frac{q^{1/2}}{\pi (s,z)}=\Delta _{s} 
\frac{1}{x(z)-x(s)},
\end{equation*}
where 
\begin{align}
\pi (s,z)& =\left[ x\left( s+1\right) -x\left( z\right) \right] \left[
x\left( s\right) -x\left( z\right) \right]  \notag \\
& =q\left[ x\left( s\right) -x\left( z-1\right) \right] \left[ x\left(
s\right) -x\left( z\right) \right] ,  \label{poly-pi}
\end{align}
is a polynomial of degree two both in the variable $s$ and $z$.

Hence, 
\begin{equation*}
\nabla S_{q}(z)=-\left\langle \mathcal{U},\Delta _{s}\frac{1}{x(z)-x(s)}
\right\rangle ,
\end{equation*}
which implies that equation (\ref{ext-def-2}) is fulfilled.

The following theorem establishes the main result in this paper.

\begin{theorem}
\label{main-th} The polynomial sequence $\{P_{n}(s)\}_{n\geq 0}$ orthogonal
with respect to the moment functional $\mathcal{U}$ is classical if and only
if the Stieltjes function (\ref{q-stieltjes}) (in terms of $q$-falling
factorials) satisfies the following first order non-homogeneous difference
equation 
\begin{equation}
\nabla \left[ (\sigma (s)+\tau (s)\bigtriangledown x(s+1/2))S_{q}(s)\right] =\tau (s)S_{q}(s)+C_{q},\quad C_{q}\in 
\mathbb{R}\backslash \left\{ 0\right\},  \label{main-eq-paper}
\end{equation}
where the constant 
\begin{equation*}
C_{q}=\left(a_{2}q^{-1/2}+\frac{1}{2}b_{1}q^{-1}(q-1)-b_{1}\right)u_{0}^{q},
\end{equation*}
depends on the polynomial coefficients $a_2$ and $b_1$ of $\widetilde{\sigma 
}(s)$ and $\widetilde{\tau}(s)$, respectively (see formulas (\ref
{sigma-coeff})-(\ref{tau-coeff})).
\end{theorem}

\emph{Proof}. Let us introduce the linear functional $\mathcal{V}=p(s) 
\mathcal{U}$, where 
\begin{align*}
p\left( s\right) & =\sigma (s)+\tau \left( s\right) \bigtriangledown x(s+ 
\mbox{\scriptsize $\frac{1}{2}$})  \notag \\
& =\alpha x^{2}\left( s\right) +\beta x\left( s\right) +\gamma ,
\label{poly-p}
\end{align*}
and 
\begin{align*}
\alpha & =a_{2}+\frac{1}{2}b_{1}q^{-1/2}\left( q-1\right) , \\
\beta & =a_{1}+\frac{1}{2}b_{1}q^{-1/2}+\frac{1}{2}b_{0}q^{-1/2}\left(
q-1\right) , \\
\gamma & =a_{0}+\frac{1}{2}b_{0}q^{-1/2}.
\end{align*}
Here we have used that 
\begin{align*}
\bigtriangledown x\left( s+\mbox{\scriptsize $\frac{1}{2}$}\right) =&
q^{-1/2}\left( q-1\right) \left[ x\left( s\right) +\frac{1}{q-1}\right] \\
=& q^{-1/2}\left( q-1\right) x\left( s\right) +q^{-1/2}.
\end{align*}
Notice that 
\begin{equation*}
\frac{x\left( s\right) -x\left( s-1\right) }{x\left( s+1/2\right) -x\left(
s-1/2\right) }=q^{-1/2}.
\end{equation*}
Hence, 
\begin{align}
\nabla p\left( s\right) & =\alpha \nabla x^{2}\left( s\right) +\beta \nabla
x\left( s\right)  \notag \\
& =\alpha \left[ \frac{x^{2}\left( s\right) -x^{2}\left( s-1\right) }{
x\left( s+1/2\right) -x\left( s-1/2\right) }\right] +\beta \left[ \frac{
x\left( s\right) -x\left( s-1\right) }{x\left( s+1/2\right) -x\left(
s-1/2\right) }\right]  \notag \\
& =\alpha q^{-1/2}\left[ x\left( s\right) +x\left( s-1\right) \right] +\beta
q^{-1/2}.  \label{delta-preferece1}
\end{align}
Thus, 
\begin{equation*}
\nabla p\left( s\right) =\alpha q^{-1/2}\left[ x\left( s\right) +x\left(
s-1\right) \right] +\beta q^{-1/2}.
\end{equation*}
Now, if one assumes that $\mathcal{U}$ is classic (see equation (\ref
{classic-func})), or equivalently, 
\begin{equation}
\nabla \left[ p\left( s\right) \mathcal{U}\right] =\tau (s)\mathcal{U},
\label{implication-backward}
\end{equation}
one gets 
\begin{equation}
\left\langle \nabla \left[ p\left( s\right) \mathcal{U}\right] ,\frac{1}{
x\left( z\right) -x\left( s\right) }\right\rangle =\left\langle \tau (s) 
\mathcal{U},\frac{1}{x\left( z\right) -x\left( s\right) }\right\rangle .
\label{ref2}
\end{equation}
According to Definition \ref{def-prod-pU}, one gets for the right hand-side
of the equation (\ref{ref2}) the following relation 
\begin{align}
\left\langle \tau (s)\mathcal{U},\frac{1}{x\left( z\right) -x\left( s\right) 
}\right\rangle & =\left\langle \mathcal{U},\frac{\tau (s)}{x\left( z\right)
-x\left( s\right) }\right\rangle =\left\langle \mathcal{U},\frac{\tau
(s)-\tau (z)+\tau (z)}{x\left( z\right) -x\left( s\right) }\right\rangle 
\notag \\
& =\tau (z)S_{q}\left( z\right) -b_{1}\left\langle \mathcal{U}
,1\right\rangle =\tau (z)S_{q}\left( z\right) -b_{1}u_{0}^{q}.  \label{ref3}
\end{align}
By Lemma \ref{diff-1-Stieltjes} and taking into account the formula before (\ref{poly-pi}), we have 
\begin{equation*}
\left\langle \nabla _{s}\mathcal{V},\frac{1}{x\left( z\right) -x\left(
s\right) }\right\rangle =-\left\langle \mathcal{V},\frac{q^{1/2}}{\pi (s,z)}
\right\rangle .
\end{equation*}
Therefore, from Definition \ref{def-prod-pU}, for the left-hand side of
equation (\ref{ref2}) one gets 
\begin{equation*}
\left\langle \nabla _{s}\left[ p\left( s\right) \mathcal{U}\right] ,\frac{1}{
x\left( z\right) -x\left( s\right) }\right\rangle =-q^{1/2}\left\langle 
\mathcal{U},\frac{p\left( s\right) }{\pi (s,z)}\right\rangle .
\end{equation*}
Thus, 
\begin{align*}
\left\langle \mathcal{U},\frac{p\left( s\right) }{\pi (s,z)}\right\rangle &
=\left\langle \mathcal{U},\frac{p\left( s\right) -p\left( z-1\right)
+p\left( z-1\right) }{\pi (s,z)}\right\rangle \\
& =-q^{-1/2}p\left( z-1\right) \nabla S_{q}\left( z\right) +\left\langle 
\mathcal{U},\frac{p\left( s\right) -p\left( z-1\right) }{\pi (s,z)}
\right\rangle .
\end{align*}
Since, 
\begin{equation*}
p\left( s\right) -p\left( z-1\right) =\alpha \left[ x^{2}\left( s\right)
-x^{2}\left( z-1\right) \right] +\beta \left[ x\left( s\right) -x\left(
z-1\right) \right] .
\end{equation*}
From expression (\ref{poly-pi}) we have 
\begin{equation*}
\frac{p\left( s\right) -p\left( z-1\right) }{\pi (s,z)}=\frac{\alpha }{q} 
\frac{x(s)+x(z-1)}{x(s)-x(z)}+\frac{\beta }{q}\frac{1}{x(s)-x(z)}.
\end{equation*}
Thus, 
\begin{align*}
\left\langle \mathcal{U},\frac{p\left( s\right) -p\left( z-1\right) }{\pi
(s,z)}\right\rangle & =\frac{\alpha }{q}\left\langle \mathcal{\ U},\frac{
x\left( s\right) +x\left( z-1\right) }{x\left( s\right) -x\left( z\right) }
\right\rangle -\frac{\beta }{q}\left\langle \mathcal{U},\frac{1}{x\left(
z\right) -x\left( s\right) }\right\rangle \\
& =\frac{\alpha }{q}u_{0}^{q}-\left( \frac{\alpha }{q}\left[ x\left(
z\right) +x\left( z-1\right) \right] +\frac{\beta }{q}\right) S_{q}\left(
z\right) .
\end{align*}
Taking into account (\ref{delta-preferece1}) the above expression transforms
into the equation 
\begin{equation*}
\left\langle \mathcal{U},\frac{p\left( s\right) -p\left( z-1\right) }{\pi
(s,z)}\right\rangle =\alpha q^{-1}u_{0}^{q}-q^{-1/2}S_{q}\left( z\right)
\nabla p\left( z\right) .
\end{equation*}
Therefore, the following expression for the left-hand side of (\ref{ref2}) 
\begin{align}
\left\langle \nabla \left[ p\left( s\right) \mathcal{U}\right] ,\frac{1}{
x\left( z\right) -x\left( s\right) }\right\rangle & =p\left( z-1\right)
\nabla S_{q}\left( z\right) +S_{q}\left( z\right) \nabla p\left( z\right)
-\alpha q^{-1/2}u_{0}^{q}  \notag \\
& =\nabla \left[ p\left( z\right) S_{q}\left( z\right) \right] -\alpha
q^{-1/2}u_{0}^{q}.  \label{ref4}
\end{align}
holds.

Then, by using (\ref{ref2}), (\ref{ref3}), and (\ref{ref4}) one gets 
\begin{eqnarray*}
\nabla \left[ p\left( z\right) S_{q}\left( z\right) \right] &=&\tau
(z)S_{q}\left( z\right) +\left( \alpha q^{-1/2}-b_{1}\right) u_{0}^{q} \\
&=&\tau (z)S_{q}\left( z\right) +\left[ a_{2}q^{-1/2}+\frac{1}{2}
b_{1}q^{-1}\left( q-1\right) -b_{1}\right] u_{0}^{q},
\end{eqnarray*}
and equation (\ref{main-eq-paper}) holds.

Notice that the statement holds since the converse implication, i.e. (\ref
{main-eq-paper}) $\Longrightarrow $ (\ref{implication-backward}), follows by
the above chain of equalities but proceeding in reverse order.

\begin{remark} Theorem \ref{main-th} can be proved in a more elementary way, i.e., without using the assertions of Lemma \ref{lemma-first} and Lemma \ref{diff-1-Stieltjes}, respectively. However, this procedure requires
some cumbersome calculations.
\end{remark}
Below we highlight the main ideas of this procedure. Assume that $\mathcal{U}$ verifies (\ref{classic-func}), or equivalently
\begin{equation*}
\left\langle \nabla \left[ (\sigma (s)+
\tau (s)\bigtriangledown x(s+1/2))
\mathcal{U}\right] ,\left[ s\right] _{q}^{(k)}\right\rangle =\left\langle
\tau (s)\mathcal{U},\left[ s\right] _{q}^{(k)}\right\rangle.
\end{equation*}
Hence,
\begin{equation}
\left\langle \mathcal{U},\left[ \sigma (s)+\tau (s)\bigtriangledown x(s+1/2)\right]
q^{1-k/2}\left[ k\right] _{q}\left[ s\right] _{q}^{(k-1)}+\tau (s)\left[s
\right] _{q}^{(k)}\right\rangle =0.  \label{refl}
\end{equation}
From equation (\ref{refl}), after some straightforward computations, one gets a three-term recurrence relation involved the following $q$-moments $u_{k-1}^{q}$, $u_{k}^{q}$, and $u_{k+1}^{q}$. The resulting equation can be expressed in the following convenient way
\begin{equation}
\Gamma _{k}+\Psi _{k}=\Xi _{k},  \label{refp}
\end{equation}
where
\begin{align*}
&\Gamma _{k}=-\left[ a_{2}q^{\frac{3k}{2}-1}+\frac{1}{2}b_{1}(q-1)q^{\frac{
3k-3}{2}}\right] \left[ k+2\right] _{q}u_{k+1}^{q} -\left[k+1\right] _{q}\left[2a_{2}q^{k-1}\left[k\right]_{q}\right.\\
&\left.-a_{2}q^{\frac{3k-3}{2}}+a_{1}q^{\frac{k-1}{2}}
+\frac{1}{2}b_{1}q^{\frac{k}{2}-1}(q^{k-1}+q^{k}-1)
+\frac{1}{2}b_{0}q^{\frac{k}{2}-1}(q-1)\right] u_{k}^{q}\\
&-\left[ k\right] _{q}\left[ a_{2}q^{\frac{k}{2}-1}\left[ k-1\right]
_{q}^{2}+(a_{1}+\frac{1}{2}b_{1}q^{k-\frac{3}{2}})\left[ k-1\right]
_{q}+(a_{0}+\frac{1}{2}b_{0}q^{k-\frac{3}{2}})q^{1-\frac{k}{2}}\right]
u_{k-1}^{q},
\end{align*}
\begin{align*}
\Psi _{k}&=\left[ a_{2}q^{k-\frac{3}{2}}(q+1)+\frac{1}{2}
b_{1}q^{k-2}(q^{2}-1)\right] u_{k+1}^{q}+\left\{ a_{2}\left[ q^{\frac{k}{2}-2}(q+1)\left[ k\right] _{q}
-q^{-\frac{3}{2}}\right]\right.\\
&\left.+a_{1}q^{-\frac{1}{2}}+\frac{1}{2}b_{1}\left[ q^{\frac{k-1}{2}}(1-q^{-2})
\left[ k\right] _{q}+q^{-2}\right] +\frac{1}{2}\frac{(q-1)}{q}b_{0}
\right\} u_{k}^{q},
\end{align*}
and
\begin{equation*}
\Xi _{k}=b_{1}q^{k}u_{k+1}^{q}+(b_{1}q^{\frac{k-1}{2}}\left[ k\right]
_{q}+b_{0})u_{k}^{q}.
\end{equation*}
Dividing $\Gamma _{k}$, $\Psi _{k}$, and $\Xi _{k}$ by $q^{k}\left[ s\right] _{q}^{(k+1)}$ and summing from $k=0$ to $\infty$, after tedious calculations, one obtains
\begin{align}
\sum\limits_{k\geq 0}\frac{\Gamma _{k}}{q^{k}\left[ s\right] _{q}^{(k+1)}}&=
\left[ a_{2}q^{-\frac{3}{2}}+\frac{1}{2}b_{1}q^{-2}(q-1)\right] u_{0}^{q}\nonumber\\
&+\left[
\sigma (s-1)+\tau (s-1)\bigtriangledown x(s-1/2)\right] \nabla S_{q}(s),
\label{ref5}
\end{align}
\begin{align}
\sum\limits_{k\geq 0}\frac{\Psi _{k}}{q^{k}\left[ s\right] _{q}^{(k+1)}}&=-
\left[ a_{2}q^{-\frac{3}{2}}(q+1)+\frac{1}{2}b_{1}(1-q^{-2})\right] u_{0}^{q}\nonumber\\
&+\nabla \left[ \sigma (s)+\tau (s)\bigtriangledown x(s+1/2)\right] S_{q}(s),
\label{ref6}
\end{align}
and 
\begin{equation}
\sum\limits_{k\geq 0}\frac{\Xi _{k}}{q^{k}\left[ s\right] _{q}^{(k+1)}}
=-b_{1}u_{0}^{q}+\tau (s)S_{q}(s).  \label{ref7}
\end{equation}%
Recall that from (\ref{refp}) one has
\begin{equation*}
\sum\limits_{k\geq 0}\frac{\Gamma _{k}}{q^{k}\left[ s\right] _{q}^{(k+1)}}
+\sum\limits_{k\geq 0}\frac{\Psi _{k}}{q^{k}\left[ s\right] _{q}^{(k+1)}}
=\sum\limits_{k\geq 0}\frac{\Xi _{k}}{q^{k}\left[ s\right] _{q}^{(k+1)}}.
\end{equation*}
Then, relations (\ref{ref5})-(\ref{ref7}) yields the equation (\ref{main-eq-paper}). The converse implication follows by the chain of obtained equalities but proceeding in reverse order.

Observe that the more standard properties that characterize the $q$-classical orthogonal polynomials can be derived from a distributional
difference equation (\ref{implication-backward}). For this observation we
refer to the approach given in \cite{me-al-ma} in the framework of a pure
algebraic approach --used also in the above proof. Accordingly, from Theorem 
\ref{main-th} follows that the $q$-classical orthogonal polynomials verify
the hy\-per\-geo\-me\-tric-type difference equation (\ref{eqdif}), the
three-term recurrence relation (\ref{ttrrc}); their finite differences
constitute an orthogonal polynomial family (\ref{ortho-diffs}) and they can
be expressed by a Rodrigues-type formula (\ref{rodrigues-t}), among other
aforementioned properties.

\section{Examples\label{examples}}

\noindent In this Section, we will give an explicit expression -case by
case- for the moment sequence (\ref{q-moment-u}) as well as for the solution
of (\ref{main-eq-paper}) in terms of $q$-hypergeometric series, i.e. of
expression (\ref{q-stieltjes}).

\subsection{$q$-Charlier case}

The $q$-Charlier polynomials are orthogonal with respect to the $q$-moment
functional (\ref{def-q-functional}) defined by the weight function (see e.g. 
\cite[eq. (87)]{Arv1}) 
\begin{equation*}
\rho (s)=\frac{\mu ^{s}}{e_{q}\left[ (1-q)\mu \right] \Gamma _{q}(s+1)}
,\quad \mu >0,\quad 0<(1-q)\mu <1,
\end{equation*}
where $s\in \left[ 0,\infty \right) $, and $e_{q}(z)$ denotes the $q$
-analogue of the exponential function (for details, see \cite{Gasper}).

Observe that the above function $\rho (s)$ is a solution of equation (\ref
{pearson2}) for the polynomial coefficients: $\tau (s)=\mu
q^{3/2}-q^{1/2}x(s)$, $\sigma (s)=q^{s}x(s)$, and 
\begin{equation*}
\sigma (s)+\tau \left( s\right) \bigtriangledown x(s+ 
\mbox{\scriptsize
$\frac{1}{2}$})=\mu q^{s+1}.
\end{equation*}
These polynomial coefficients (see e.g. \cite[eq. (86)]{Arv1}) are involved
in the characteristic equation (\ref{main-eq-paper}) as well as the moment $
u_{0}^{q}$. For determining this moment we now compute explicitly all the
moments of (\ref{def-q-functional}) for the above choice of $\rho (s)$.
According to (\ref{q-moment-u}) one has 
\begin{eqnarray*}
u_{k}^{q} &=&\frac{q^{-1/2}}{e_{q}\left[ (1-q)\mu \right] }\sum_{s\geq k} 
\left[ s\right] _{q}^{(k)}\frac{\widetilde{\mu }^{s}}{\Gamma _{q}(s+1)}= 
\frac{q^{-1/2}}{e_{q}\left[ (1-q)\mu \right] }\sum_{s\geq k}\frac{\widetilde{
\mu }^{s}}{\displaystyle\frac{(q;q)_{s-k}}{(1-q)^{s-k}}} \\
&=&\frac{q^{-1/2}\widetilde{\mu }^{k}}{e_{q}\left[ (1-q)\mu \right] }
\sum_{n\geq 0}\frac{\widetilde{\mu }^{n}}{\displaystyle\frac{(q;q)_{n}}{
(1-q)^{n}}},\quad \mbox{where}\quad \widetilde{\mu }=q\mu .
\end{eqnarray*}
Thus, the explicit expression for the $k$-th moment associated with (\ref
{def-q-functional}) is as follows 
\begin{equation*}
u_{k}^{q}=\frac{\widetilde{\mu }^{k}e_{q}\left[ (1-q)\widetilde{\mu }\right] 
}{q^{1/2}e_{q}\left[ (1-q)\mu \right] },\quad k=0,1,\dots
\end{equation*}
In particular 
\begin{equation*}
u_{0}^{q}=\frac{e_{q}\left[ (1-q)\widetilde{\mu }\right] }{q^{1/2}e_{q}\left[
(1-q)\mu \right] }.
\end{equation*}
Finally, from (\ref{q-stieltjes}) one gets 
\begin{equation*}
S_{q}(z)=\frac{u_{0}^{q}}{x\left( z\right) }\sum_{k\geq 0}\frac{\left(
q;q\right) _{k}\left( -1\right) ^{k}q^{\binom{k}{2}}}{\left(
q^{1-z};q\right) _{k}}\frac{\left[ \mu (1-q)q^{1-z}\right] ^{k}}{\left(
q;q\right) _{k}},
\end{equation*}
since the $q$-falling factorial can be rewritten as follows (see (\ref
{ref-q-falling})) 
\begin{equation}
\left[ z\right] _{q}^{\left( k+1\right) }=\left( -1\right) ^{k}\frac{x\left(
z\right) }{\left( 1-q\right) ^{k}}\left( q^{1-z};q\right) _{k}q^{\left(
z-1\right) k-\binom{k}{2}}.  \label{q-fallingT}
\end{equation}
Hence, the $q$-Stieltjes function associated with the $q$-Charlier moment
functional that solves equation (\ref{main-eq-paper}) is given in terms of
the following hypergeometric series 
\begin{equation}
S_{q}(z)=\frac{u_{0}^{q}}{x\left( z\right) }\,_{1}\varphi _{1}\left( 
\begin{array}{c|c}
q &  \\ 
& q;\mu (1-q)q^{1-z} \\ 
q^{1-z} & 
\end{array}
\right) .  \label{q-StieltjesFunction-Charlier1}
\end{equation}
On the other hand, we know that (\ref{q-StieltjesFunction-Charlier1}) must
be equal to (\ref{q-stieltjes-Suslov}). This fact is easily established if
one rewrites $\rho (s)$ as 
\begin{equation*}
\rho (s)=C_{q}\frac{\left[ (1-q)\mu \right] ^{s}}{\left( q;q\right) _{s}}
,\quad \mbox{where}\quad C_{q}=\frac{1}{e_{q}\left[ (1-q)\mu \right] },
\end{equation*}
and (\ref{q-stieltjes-Suslov}) as 
\begin{eqnarray*}
S_{q}(z) &=&\frac{C_{q}q^{-1/2}}{x\left( z\right) }\sum_{s\geq 0}\frac{
\left( q^{-z};q\right) _{s}}{\left( q^{1-z};q\right) _{s}}\frac{\left[
(1-q)\mu q\right] ^{s}}{\left( q;q\right) _{s}} \\
&=&\frac{C_{q}q^{-1/2}}{x\left( z\right) }\sum_{s\geq 0}\frac{\left(
q^{-z};q\right) _{s}\left( 0;q\right) _{s}}{\left( q^{1-z};q\right) _{s}} 
\frac{\left[ (1-q)\mu q\right] ^{s}}{\left( q;q\right) _{s}}.
\end{eqnarray*}
Here we have used that $\Gamma _{q}\left( s+1\right) =\frac{\left(
q;q\right) _{s}}{(1-q)^{s}}$.

Thus, the $q$-Stieltjes function associated with the $q$-Charlier moment
functional also has the form 
\begin{equation*}
S_{q}(z)=\frac{C_{q}q^{-1/2}}{x\left( z\right) }\,_{2}\varphi _{1}\left( 
\begin{array}{c|c}
q^{-z},0 &  \\ 
& q;(1-q)\mu q \\ 
q^{1-z} & 
\end{array}
\right) .  \label{q-StieltjesFunction-Charlier2}
\end{equation*}
Now, using the Heine's transformation formula (for details, see \cite
{koekoek} page 16) 
\begin{equation*}
_{2}\varphi _{1}\left( 
\begin{array}{c|c}
a,0 &  \\ 
& q;z \\ 
c & 
\end{array}
\right) =e_{q}\left[ z\right] \,_{1}\varphi _{1}\left( 
\begin{array}{c|c}
a^{-1}c &  \\ 
& q;az \\ 
c & 
\end{array}
\right) ,\quad a\neq 0,\;\left\vert z\right\vert <1,
\end{equation*}
under the assumption $\left\vert (1-q)\mu q\right\vert <1$ one gets the
desired equality: 
\begin{eqnarray*}
S_{q}(z) &=&\frac{C_{q}q^{-1/2}}{x\left( z\right) }\,_{2}\varphi _{1}\left( 
\begin{array}{c|c}
q^{-z},0 &  \\ 
& q;(1-q)\mu q \\ 
q^{1-z} & 
\end{array}
\right) \\
&=&\frac{u_{0}^{q}}{x\left( z\right) }\,_{1}\varphi _{1}\left( 
\begin{array}{c|c}
q &  \\ 
& q;\mu (1-q)q^{1-z} \\ 
q^{1-z} & 
\end{array}
\right) .
\end{eqnarray*}

\subsection{$q$-Kravchuk case}

The $q$-Kravchuk polynomials are orthogonal with respect to the $q$-moment
functional (\ref{def-q-functional}), where the orthogonalizing weight
function 
\begin{equation*}
\rho (s)=q^{\binom{s}{2}}\frac{[N]_{q}!}{\Gamma _{q}(s+1)\Gamma _{q}(N-s+1)}
p^{s}(1-p)^{N-s},\quad 0<p<1,\quad N\in \mathbb{N},
\end{equation*}
is a solution of Pearson-type equation (\ref{pearson2}) for the choice (see \cite{arvesu-tesis} p.89)
\begin{align*}
\tau (s)& =\displaystyle\frac{q^{1/2}pq(q^{N}-1)}{1-p}-\frac{q^{1/2}\left(
p(q-1)+1\right) }{1-p}x(s) \\
\sigma (s)& =(q-1)x^{2}(s)+x(s).
\end{align*}
By using the expression $\Gamma _{q}(s+1)=q^{\binom{s}{2}/2}[s]_{q}!$, the
above function $\rho (s)$ can be rewritten as 
\begin{equation*}
\rho (s)=q^{s(N-1)/2-\binom{N}{2}/2}\frac{[N]_{q}!}{[s]_{q}![N-s]_{q}!}
p^{s}(1-p)^{N-s}.
\end{equation*}
Hence, for the moments (\ref{q-moment-u}) one gets 
\begin{eqnarray*}
u_{k}^{q} &=&q^{-\frac{1}{2}\left[ \binom{N}{2}+1\right] }
\sum_{s=k}^{N}q^{s(N+1)/2}\left[ s\right] _{q}^{(k)}\frac{[N]_{q}!}{
[s]_{q}![N-s]_{q}!}p^{s}(1-p)^{N-s} \\
&=&q^{-\frac{1}{2}\left[ \binom{k+1}{2}+\binom{N}{2}+1\right]
}\sum_{s=k}^{N}q^{s(N+k+1)/2}\frac{[N]_{q}!}{[s-k]_{q}![N-s]_{q}!}
p^{s}(1-p)^{N-s} \\
&=&q^{\frac{1}{2}\left[ 2\binom{k+1}{2}-\binom{N}{2}-1\right]
}[N]_{q}^{(k)}p^{k}\sum_{n=0}^{N-k}\frac{q^{n(N+k+1)/2}[N-k]_{q}!}{
[n]_{q}![N-k-n]_{q}!}p^{n}(1-p)^{N-n-k},
\end{eqnarray*}
where we have used the following relations 
\begin{align*}
\frac{\left[ s\right] _{q}^{(k)}}{[s]_{q}!}& =\frac{q^{ks/2-\binom{k+1}{2}
/2} }{[s-k]_{q}!}, \\
\lbrack N]_{q}!& =q^{\frac{1}{2}\left[ \binom{N-k}{2}-\binom{N}{2}\right]
}[N]_{q}^{(k)}[N-k]_{q}!.
\end{align*}
Equivalently, 
\begin{equation*}
u_{k}^{q}=q^{\frac{1}{2}\left[ 2\binom{k+1}{2}-\binom{N}{2}-1\right]
}[N]_{q}^{(k)}p^{k}\sum_{n=0}^{N-k}q^{n(2k+n+1)/2}\left[ 
\begin{array}{c}
N-k \\ 
n
\end{array}
\right] _{q}p^{n}(1-p)^{N-n-k},
\end{equation*}
where the $q$-binomial symbol is defined as (see e.g. \cite[p. 24]{Gasper}) 
\begin{equation*}
\left[ 
\begin{array}{c}
n \\ 
j
\end{array}
\right] _{q}=\frac{\Gamma _{q}\left( n+1\right) }{\Gamma _{q}\left(
k+1\right) \Gamma _{q}\left( n-k+1\right) }=q^{k\left( n-k\right) /2}\frac{
[n]_{q}!}{[k]_{q}![n-k]_{q}!}.
\end{equation*}
Using the well-known relation (see \cite[pp. 3--6]{DavidM}) 
\begin{equation*}
\sum_{j=0}^{n}q^{\binom{j+1}{2}}\left[ 
\begin{array}{c}
n \\ 
j
\end{array}
\right] _{q}x^{n-j}y^{j}=\prod_{j=0}^{n-1}\left( x+yq^{j+1}\right) ,
\end{equation*}
the above expression for the $k$-th moment transforms into 
\begin{equation*}
u_{k}^{q}=q^{\frac{1}{2}\left[ 2\binom{k+1}{2}-\binom{N}{2}-1\right]
}[N]_{q}^{(k)}p^{k}\prod_{j=0}^{N-k-1}\left( 1-p+pq^{j+k+1}\right) .
\end{equation*}
Finally 
\begin{equation*}
u_{k}^{q}=\frac{u_{0}^{q}}{\left( 1-q\right) ^{k}}\frac{\left(
q^{-N};q\right) _{k}}{\left( \displaystyle\frac{pq}{p-1};q\right) _{k}}
\left( \frac{pq^{N+1}}{p-1}\right) ^{k},
\end{equation*}
where 
\begin{equation*}
u_{0}^{q}=\frac{\left( 1-p\right) ^{N}}{\sqrt{q^{\binom{N}{2}+1}}}\left( 
\frac{pq}{p-1};q\right) _{N}.
\end{equation*}
Therefore, from (\ref{q-stieltjes}) and (\ref{q-fallingT}) one gets the
following explicit expression for the $q$-Stieltjes function associated with
the $q$-Kravchuk moment functional 
\begin{equation*}
S_{q}(z)=\frac{u_{0}^{q}}{x\left( z\right) }\sum_{k\geq 0}\frac{\left(
q^{-N};q\right) _{k}\left( q;q\right) _{k}\left( -1\right) ^{k}q^{\binom{k}{
2 }}}{\left( q^{1-z};q\right) _{k}\left( \displaystyle\frac{pq}{p-1}
;q\right) _{k}}\frac{\left( \frac{pq^{N+1-z}}{p-1}\right) ^{k}}{\left(
q;q\right) _{k}} .
\end{equation*}
This expression can be written in terms of hypergeometric series as follows 
\begin{equation}
S_{q}(z)=\frac{u_{0}^{q}}{x\left( z\right) }\,_{2}\varphi _{2}\left( 
\begin{array}{c|c}
q^{-N},q &  \\ 
& q;\displaystyle\frac{pq^{N+1-z}}{p-1} \\ 
q^{1-z},\displaystyle\frac{pq}{p-1} & 
\end{array}
\right) .  \label{q-Stieltjes functionK1}
\end{equation}
Previously, we have established a relationship between (\ref{q-Stieltjes
functionK1}) and (\ref{q-stieltjes-Suslov}). Now, aimed to check a similar
relationship between (\ref{q-Stieltjes functionK1}) and (\ref
{q-stieltjes-Suslov}) we rewrite $\rho (s)$ as 
\begin{equation}
\rho (s)=C_{q}\frac{\left( q^{-N};q\right) _{s}}{\left( q;q\right) _{s}}
\left( \frac{pq^{N}}{p-1}\right) ^{s},  \label{kravchuk-rho-2}
\end{equation}
where 
\begin{equation*}
C_{q}=\left( 1-p\right) ^{N}\frac{\left[ N\right] _{q}!}{\Gamma _{q}\left(
N+1\right) }=\left( 1-p\right) ^{N}q^{-\binom{N}{2}/2}.
\end{equation*}
In (\ref{kravchuk-rho-2}) we have taken into account the following relation 
\begin{equation*}
\frac{\Gamma _{q}(N-s+1)}{\Gamma _{q}(N+1)}=\frac{\left( -1\right)
^{s}\left( 1-q\right) ^{s}q^{\binom{s}{2}-Ns}}{\left( q^{-N};q\right) _{s}}.
\end{equation*}
Thus, from (\ref{q-stieltjes-Suslov}) one gets 
\begin{equation*}
S_{q}(z)=\frac{C_{q}q^{-1/2}}{x\left( z\right) }\sum_{s\geq 0}\frac{\left(
q^{-N};q\right) _{s}\left( q^{-z};q\right) _{s}}{\left( q^{1-z};q\right)
_{s} }\frac{\left( \frac{pq^{N+1}}{p-1}\right) ^{s}}{\left( q;q\right) _{s}}.
\end{equation*}
Accordingly, the $q$-Stieltjes function associated with the $q$-Kravchuk
moment functional also has the form 
\begin{equation*}
S_{q}(z)=\frac{C_{q}q^{-1/2}}{x\left( z\right) }\,_{2}\varphi _{1}\left( 
\begin{array}{c|c}
q^{-N},q^{-z} &  \\ 
& q;\displaystyle\frac{pq^{N+1}}{p-1} \\ 
q^{1-z} & 
\end{array}
\right) .
\end{equation*}
Now, using the Jackson's transformation formula (for details, see \cite
{koekoek} p. 15) 
\begin{equation}
_{2}\varphi _{1}\left( 
\begin{array}{c|c}
a,b &  \\ 
& q;z \\ 
c & 
\end{array}
\right) =\frac{\left( az;q\right) _{\infty }}{\left( z;q\right) _{\infty }}
\,_{2}\varphi _{2}\left( 
\begin{array}{c|c}
a,b^{-1}c &  \\ 
& q;bz \\ 
c,az & 
\end{array}
\right) ,\quad b\neq 0,  \label{Jackson'sTransf}
\end{equation}
one obtains the aforementioned relationship 
\begin{eqnarray*}
S_{q}(z) &=&\frac{C_{q}q^{-1/2}}{x\left( z\right) }\,_{2}\varphi _{1}\left( 
\begin{array}{c|c}
q^{-N},q^{-z} &  \\ 
& q;\displaystyle\frac{pq^{N+1}}{p-1} \\ 
q^{1-z} & 
\end{array}
\right) \\
&=&\frac{C_{q}q^{-1/2}}{x\left( z\right) }\frac{\left( \frac{pq}{p-1}
;q\right) _{\infty }}{\left( \frac{pq}{p-1}q^{N};q\right) _{\infty }}
\,_{2}\varphi _{2}\left( 
\begin{array}{c|c}
q^{-N},q &  \\ 
& q;\displaystyle\frac{pq^{N+1-z}}{p-1} \\ 
q^{1-z},\displaystyle\frac{pq}{p-1} & 
\end{array}
\right) \\
&=&\frac{u_{0}^{q}}{x\left( z\right) }\,_{2}\varphi _{2}\left( 
\begin{array}{c|c}
q^{-N},q &  \\ 
& q;\displaystyle\frac{pq^{N+1-z}}{p-1} \\ 
q^{1-z},\displaystyle\frac{pq}{p-1} & 
\end{array}
\right) .
\end{eqnarray*}

\subsection{$q$-Meixner case}

The $q$-Meixner polynomials are orthogonal with respect to the $q$-moment
functional (\ref{def-q-functional}) defined by the weight function (see e.g. 
\cite[p. 314]{Arv1}) 
\begin{equation}
\rho (s)=\frac{\mu ^{s}\Gamma _{q}(\gamma +s)}{\Gamma _{q}(\gamma )\Gamma
_{q}(s+1)}=\frac{\mu ^{s}}{q^{\binom{s}{2}/2}[s]_{q}!}\frac{(q^{\gamma
};q)_{s}}{(1-q)^{s}},\quad 0<\mu <1,\quad \gamma >0.  \label{rho-m}
\end{equation}
Recall that $\rho (s)$ is a solution of Pearson-type difference equation for
the polynomial coefficients (see e.g. \cite[eq. (68)]{Arv1}, up to the
factor $q-1$): 
\begin{align*}
\tau (s)& =q^{1/2}(\mu q^{\gamma +1}-1)x(s)+\mu q^{\frac{\gamma +2}{2}
}[\gamma ]_{q}, \\
\sigma (s)& =(q-1)x(s)^{2}+x(s).
\end{align*}
For the above orthogonalizing weight function and (\ref{def-q-functional}),
let us compute the associated $q$-moments (\ref{q-moment-u}) 
\begin{eqnarray*}
u_{k}^{q} &=&q^{-1/2}\sum_{s\geq 0}q^{s-\binom{s}{2}/2}\left[ s\right]
_{q}^{(k)}\frac{\mu ^{s}}{[s]_{q}!}\frac{(q^{\gamma };q)_{s}}{(1-q)^{s}} \\
&=&q^{-\frac{1}{2}\left[ \binom{k+1}{2}+1\right] }\sum_{s\geq k}q^{s(k+2)/2- 
\binom{s}{2}/2}\frac{\mu ^{s}}{[s-k]_{q}!}\frac{(q^{\gamma };q)_{s}}{
(1-q)^{s}} \\
&=&q^{-\frac{1}{2}\left[ \binom{k+1}{2}+1-k\left( k+2\right) \right] }\mu
^{k}\sum_{n\geq 0}q^{n(k+2)/2-\binom{n+k}{2}/2}\frac{\mu ^{n}}{[n]_{q}!} 
\frac{(q^{\gamma };q)_{n+k}}{(1-q)^{n+k}}.
\end{eqnarray*}
Considering the following elementary relations 
\begin{align*}
n(k+2)-\binom{n+k}{2}& =n(5-n)/2-\binom{k}{2}, \\
(q^{\gamma };q)_{n+k}& =(q^{\gamma };q)_{k}(q^{\gamma +k};q)_{n},
\end{align*}
one gets 
\begin{eqnarray*}
u_{k}^{q} &=&q^{k-1/2}\frac{(q^{\gamma };q)_{k}}{(1-q)^{k}}\mu
^{k}\sum_{n\geq 0}q^{n(5-n)/4}\frac{\mu ^{n}}{[n]_{q}!}\frac{(q^{\gamma
+k};q)_{n}}{(1-q)^{n}} \\
&=&q^{k-1/2}\frac{\Gamma _{q}(\gamma +k)}{\Gamma _{q}(\gamma )}\mu
^{k}\sum_{n\geq 0}q^{n(5-n)/4}\frac{\mu ^{n}}{[n]_{q}!}\frac{(q^{\gamma
+k};q)_{n}}{(1-q)^{n}}.
\end{eqnarray*}
Since $[n]_{q}!=\displaystyle q^{-\binom{n}{2}/2}\frac{(q;q)_{n}}{(1-q)^{n}}$
, the following relation 
\begin{equation*}
u_{k}^{q}=q^{-1/2}\frac{\Gamma _{q}(\gamma +k)}{\Gamma _{q}(\gamma )} 
\widetilde{\mu }^{k}\sum_{n\geq 0}\frac{\widetilde{\mu }^{n}}{(q;q)_{n}}
(q^{\gamma +k};q)_{n},\quad \mbox{where}\quad \widetilde{\mu }=q\mu ,
\end{equation*}
holds. Thus, by the $q$-binomial theorem (see e.g. \cite[eq. 1.3.2]{Gasper}
), 
\begin{equation*}
u_{k}^{q}=q^{-1/2}\frac{\left( \mu q^{\gamma +1};q\right) _{\infty }}{\left(
\mu q;q\right) _{\infty }}\frac{\Gamma _{q}(\gamma +k)}{\Gamma _{q}(\gamma )}
\frac{\widetilde{\mu }^{k}}{\left( \mu q^{\gamma +1};q\right) _{k}}.
\end{equation*}
In particular 
\begin{equation*}
u_{0}^{q}=q^{-1/2}\frac{\left( \mu q^{\gamma +1};q\right) _{\infty }}{\left(
\mu q;q\right) _{\infty }}.
\end{equation*}
In order to obtain an expression for the corresponding $q$-Stieltjes
function associated with the $q$-Meixner moment functional one uses the
following elementary relation 
\begin{equation}
\frac{\Gamma _{q}(\gamma +k)}{\Gamma _{q}(\gamma )}=\frac{\left( q^{\gamma
};q\right) _{k}}{\left( 1-q\right) ^{k}},  \label{q-Gamma-m}
\end{equation}
as well as (\ref{q-stieltjes}) and (\ref{q-fallingT}), respectively. Thus, 
\begin{equation*}
S_{q}(z)=\frac{u_{0}^{q}}{x\left( z\right) }\sum_{k\geq 0}\frac{\left(
q^{\gamma };q\right) _{k}\left( q;q\right) _{k}\left( -1\right) ^{k}q^{ 
\binom{k}{2}}}{\left( q^{1-z};q\right) _{k}\left( \mu q^{\gamma +1};q\right)
_{k}}\frac{\left( \mu q^{1-z}\right) ^{k}}{\left( q;q\right) _{k}}.
\end{equation*}
Equivalently, in terms of hypergeometric series 
\begin{equation}
S_{q}(z)=\frac{u_{0}^{q}}{x\left( z\right) }\,_{2}\varphi _{2}\left( 
\begin{array}{c|c}
q^{\gamma },q &  \\ 
& q;\mu q^{1-z} \\ 
q^{1-z},\mu q^{\gamma +1} & 
\end{array}
\right) .  \label{q-StieltjesFunctionM1}
\end{equation}
Again, to establish an explicit relationship between (\ref
{q-StieltjesFunctionM1}) and (\ref{q-stieltjes-Suslov}) we rewrite (\ref
{rho-m}) --based on equation (\ref{q-Gamma-m})-- as 
\begin{equation*}
\rho (s)=\mu ^{s}\frac{\left( q^{\gamma };q\right) _{s}}{\left( q;q\right)
_{s}}.
\end{equation*}
Thus, from (\ref{q-stieltjes-Suslov}) one obtains the following $q$
-Stieltjes function associated with the $q$-Meixner moment functional 
\begin{equation*}
S_{q}(z)=\frac{q^{-1/2}}{x\left( z\right) }\sum_{s\geq 0}\frac{\left(
q^{\gamma };q\right) _{s}(q^{-z};q)_{s}}{(q^{1-z};q)_{s}}\frac{\left( \mu
q\right) ^{s}}{\left( q;q\right) _{s}},
\end{equation*}
or equivalently, in terms of hypergeometric series 
\begin{equation*}
S_{q}(z)=\frac{q^{-1/2}}{x\left( z\right) }\,_{2}\varphi _{1}\left( 
\begin{array}{c|c}
q^{\gamma };q^{-z} &  \\ 
& q;\mu q \\ 
q^{1-z} & 
\end{array}
\right) .
\end{equation*}
Again, using the Jackson's transformation formula (\ref{Jackson'sTransf}) we
have the following relation 
\begin{eqnarray*}
S_{q}(z) &=&\frac{q^{-1/2}}{x\left( z\right) }\,_{2}\varphi _{1}\left( 
\begin{array}{c|c}
q^{\gamma };q^{-z} &  \\ 
& q;\mu q \\ 
q^{1-z} & 
\end{array}
\right) \\
&=&\frac{q^{-1/2}}{x\left( z\right) }\frac{\left( \mu q^{\gamma +1};q\right)
_{\infty }}{\left( \mu q;q\right) _{\infty }}\,_{2}\varphi _{2}\left( 
\begin{array}{c|c}
q^{\gamma },q &  \\ 
& q;\mu q^{1-z} \\ 
q^{1-z},\mu q^{\gamma +1} & 
\end{array}
\right) \\
&=&\frac{u_{0}^{q}}{x\left( z\right) }\,_{2}\varphi _{2}\left( 
\begin{array}{c|c}
q^{\gamma },q &  \\ 
& q;\mu q^{1-z} \\ 
q^{1-z},\mu q^{\gamma +1} & 
\end{array}
\right) .
\end{eqnarray*}

\subsection{$q$-Hahn case}

The $q$-Hahn polynomials are orthogonal with respect to the $q$-moment
functional (\ref{def-q-functional}) defined by the orthogonalizing weight in
the interval $\left[ 0,N-1\right] $ (see \cite[Table 4.1]{Arv2}) 
\begin{equation*}
\rho (s)=q^{\left( \frac{\alpha +\beta }{2}\right) s}\frac{\widetilde{\Gamma 
}_{q}(s+\beta +1)\widetilde{\Gamma }_{q}(N+\alpha -s)}{\widetilde{\Gamma }
_{q}(s+1)\widetilde{\Gamma }_{q}(N-s)},\quad \alpha ,\beta >-1,\quad N\in 
\mathbb{N},
\end{equation*}
where $\tilde{\Gamma}(s)=q^{-(s-1)(s-2)/4}f(s;q)$ if $0<q<1$, or $\tilde{
\Gamma}(s)=f(s;q^{-1})$ if $q>1$, being $f(s;q)=(1-q)^{1-s}\prod_{k\geq
0}(1-q^{k+1})/\prod_{k\geq 0}(1-q^{s+k})$. The above function $\rho (s)$ is
a solution of equation (\ref{pearson2}) for the polynomial coefficients (see 
\cite{Arv2}): 
\begin{align*}
\displaystyle\sigma (s)& =-q^{-(N+\alpha )/2}x^{2}(s)+q^{-1/2}[N+\alpha
]_{q}x(s), \\
\displaystyle\tau (s)& =-q^{-(\beta+2-N )/2}[\alpha +\beta
+2]_{q}x(s)+q^{\alpha +\beta+1}[\beta +1]_{q}[N-1]_{q}.
\end{align*}
For calculating the $q$-moments (\ref{q-moment-u}) we use 
\begin{eqnarray*}
u_{k}^{q} &=&q^{-1/2}\sum_{s=k}^{N-1}q^{\left( \frac{\alpha +\beta +2}{2}
\right) s}\left[ s\right] _{q}^{(k)}\frac{\widetilde{\Gamma }_{q}(s+\beta
+1) \widetilde{\Gamma }_{q}(N+\alpha -s)}{[s]_{q}![N-1-s]_{q}!} \\
&=&q^{-\frac{1}{2}\left[ \binom{k+1}{2}+1\right] }\sum_{s=k}^{N-1}q^{\left( 
\frac{\alpha +\beta +k+2}{2}\right) s}\frac{\widetilde{\Gamma }_{q}(s+\beta
+1)\widetilde{\Gamma }_{q}(N+\alpha -s)}{[s-k]_{q}![N-1-s]_{q}!}.
\end{eqnarray*}
Hence, 
\begin{equation*}
u_{k}^{q}=q^{\chi (k,\alpha ,\beta )}\sum_{n=0}^{N-k-1}q^{n\left( \frac{
\alpha +\beta +k+2}{2}\right) }\frac{\widetilde{\Gamma }_{q}(\beta +k+1+n) 
\widetilde{\Gamma }_{q}(N+\alpha -k-n)}{[n]_{q}![N-k-1-n]_{q}!}
\end{equation*}
where $\chi (k,\alpha ,\beta )=k\left( \frac{\alpha +\beta +k+2}{2}\right) - 
\frac{1}{2}\left[ \binom{k+1}{2}+1\right] $, or equivalently 
\begin{multline*}
u_{k}^{q}=\frac{q^{\chi (k,\alpha ,\beta )}}{[N-k-1]_{q}!} \\
\sum_{n=0}^{N-k-1}q^{\psi \left( \alpha ,\beta ,n,N\right) }\left[ 
\begin{array}{c}
N-k-1 \\ 
n
\end{array}
\right] _{q}\widetilde{\Gamma }_{q}(\beta +k+1+n)\widetilde{\Gamma }
_{q}(N+\alpha -k-n).
\end{multline*}
For brevity we have introduced the notation: $\psi \left( \alpha ,\beta
,n,N,k\right) =n(\alpha +\beta +2k+3+n-N)/2$.

Taking into account 
\begin{align*}
\displaystyle\frac{\widetilde{\Gamma }_{q}(\beta +k+1+n)}{\Gamma _{q}\left(
\beta +k+1\right) }& =q^{-\left( \beta +k+n\right) \left( \beta
+k+n-1\right) /4}\frac{\left( q^{\beta +k+1};q\right) _{n}}{\left(
1-q\right) ^{n}} \\
\displaystyle\frac{\widetilde{\Gamma }_{q}\left( N-k-n+\alpha \right) }{
\Gamma _{q}\left( \alpha +1\right) }& =q^{-\left( N-k-1-n+\alpha \right)
\left( N-k-n+\alpha -2\right) /4}\frac{\left( q^{\alpha +1};q\right)
_{N-k-1-n}}{\left( 1-q\right) ^{N-k-1-n}},
\end{align*}
the $q$-moments becomes 
\begin{multline*}
u_{k}^{q}=\frac{\Gamma _{q}\left( \alpha +1\right) \Gamma _{q}\left( \beta
+k+1\right) q^{\upsilon (\alpha ,\beta ,N,k)}}{\left( 1-q\right)
^{N-k-1}[N-k-1]_{q}!} \\
\sum_{n=0}^{N-k-1}q^{n\left( \alpha +1\right) }\left[ 
\begin{array}{c}
N-k-1 \\ 
n
\end{array}
\right] _{q}\left( q^{\alpha +1};q\right) _{N-k-1-n}\left( q^{\beta
+k+1};q\right) _{n}.
\end{multline*}
where 
\begin{equation*}
2(\upsilon (\alpha ,\beta ,N,k)+1)=\alpha \left( 2k-N+1\right) +N\left(
k+1\right) -\binom{k}{2}-\binom{N}{2}-\binom{\alpha }{2}-\binom{\beta }{2}.
\end{equation*}
Using the following well-known relation (see \cite[p. 25]{Gasper}) 
\begin{equation*}
\sum_{j=0}^{n}\left[ 
\begin{array}{c}
n \\ 
j
\end{array}
\right] _{q}\left( a;q\right) _{n-j}\left( b;q\right) _{j}a^{j}=\left(
ab;q\right) _{n},
\end{equation*}
one gets 
\begin{equation*}
u_{k}^{q}=\frac{\Gamma _{q}\left( \alpha +1\right) \Gamma _{q}\left( \beta
+k+1\right) q^{\upsilon (\alpha ,\beta ,N,k)}}{\left( 1-q\right)
^{N-k-1}[N-k-1]_{q}!}\left( q^{\alpha +\beta +k+2};q\right) _{N-k-1}.
\end{equation*}
In particular 
\begin{equation*}
u_{0}^{q}=\frac{\Gamma _{q}\left( \alpha +1\right) \Gamma _{q}\left( \beta
+1\right) q^{\upsilon (\alpha ,\beta ,N,0)}}{\left( 1-q\right)
^{N-1}[N-1]_{q}!}\left( q^{\alpha +\beta +2};q\right) _{N-1}.
\end{equation*}
On the other hand, using relations 
\begin{align*}
\displaystyle& \frac{\Gamma _{q}(\beta +k+1)}{\Gamma _{q}\left( \beta
+1\right) }=\frac{\left( q^{\beta +1};q\right) _{k}}{\left( 1-q\right) ^{k}}
\\
\displaystyle& \frac{\widetilde{\Gamma }_{q}\left( N-k\right) }{\Gamma
_{q}\left( N\right) }=q^{-\left( N-k-1\right) \left( N-k-2\right) /4}\frac{
\left( -1\right) ^{k}\left( 1-q\right) ^{k}q^{\left( 1-N\right) k+\binom{k}{
2 }}}{\left( q^{1-N};q\right) _{k}}, \\
\displaystyle& \left( q^{\alpha +\beta +k+2};q\right) _{N-k-1}=\frac{\left(
q^{\alpha +\beta +2};q\right) _{N-1}}{\left( q^{\alpha +\beta +2};q\right)
_{k}}, \\
\displaystyle& \Gamma _{q}\left( N\right) =q^{\left( N-1\right) \left(
N-2\right) /4}\left[ N-1\right] _{q}!,
\end{align*}
we obtain 
\begin{equation*}
u_{k}^{q}=u_{0}^{q}q^{k\left( N+\alpha \right) -\binom{k}{2}}\left(
-1\right) ^{k}\frac{\left( q^{\beta +1};q\right) _{k}\left( q^{1-N};q\right)
_{k}}{\left( 1-q\right) ^{k}\left( q^{\alpha +\beta +2};q\right) _{k}}.
\end{equation*}
Finally, using (\ref{q-stieltjes}) and (\ref{q-fallingT}) one gets the $q$
-Stieltjes function associated with the $q$-Hahn moment functional 
\begin{equation*}
S_{q}(z)=\frac{u_{0}^{q}}{x\left( z\right) }\sum_{k\geq 0}\frac{\left(
q^{\beta +1};q\right) _{k}\left( q^{1-N};q\right) _{k}\left( q;q\right) _{k} 
}{\left( q^{1-z};q\right) _{k}\left( q^{\alpha +\beta +2};q\right) _{k}} 
\frac{q^{k\left( N+\alpha -z\right) }}{\left( q;q\right) _{k}},
\end{equation*}
or equivalently --in terms of hypergeometric series-- 
\begin{equation*}
S_{q}(z)=\frac{u_{0}^{q}}{x\left( z\right) }\,_{3}\varphi _{2}\left( 
\begin{array}{c|c}
q^{\beta +1},q^{1-N},q &  \\ 
& q;q^{N+\alpha -z} \\ 
q^{1-z},q^{\alpha +\beta +2} & 
\end{array}
\right) .
\end{equation*}
Now, aimed to find an equivalent relation for the above $q$-Stieltjes
function one uses the relations 
\begin{equation*}
\widetilde{\Gamma }_{q}\left( s+1\right) =q^{-\binom{s}{2}/2}\frac{\left(
q;q\right) _{s}}{\left( 1-q\right) ^{s}},
\end{equation*}
and 
\begin{equation*}
\frac{\widetilde{\Gamma }_{q}\left( s+\beta +1\right) }{\Gamma _{q}\left(
\beta +1\right) } =q^{-\left( s+\beta \right) \left( s+\beta -1\right) /4} 
\frac{\left( q^{\beta +1};q\right) _{s}}{\left( 1-q\right) ^{s}},
\end{equation*}
as well as
\begin{equation*}
\frac{\widetilde{\Gamma }_{q}\left( N+\alpha -s\right) }{\Gamma _{q}\left(
N+\alpha \right) }=q^{-\left( N+\alpha -s-1\right) \left( N+\alpha
-s-2\right) /4}\frac{\left( -1\right) ^{s}\left( 1-q\right) ^{s}q^{\left(
1-N-\alpha \right) s+\binom{s}{2}}}{\left( q^{1-N-\alpha };q\right) _{s}}.
\end{equation*}
Thus, we can rewrite $\rho (s)$ as follows 
\begin{equation*}
\rho (s)=q^{\tilde{\upsilon}(\alpha ,\beta ,N)}\frac{\Gamma _{q}\left( \beta
+1\right) \Gamma _{q}\left( N+\alpha \right) }{\Gamma _{q}\left( N\right) } 
\frac{\left( q^{\beta +1};q\right) _{s}\left( q^{1-N};q\right) _{s}}{\left(
q^{1-N-\alpha };q\right) _{s}\left( q;q\right) _{s}},
\end{equation*}
where $2\tilde{\upsilon}(\alpha ,\beta ,N)=\alpha \left( 1-N\right) -\binom{
\alpha }{2}-\binom{\beta }{2}$. Then, from (\ref{q-stieltjes-Suslov}) one
gets 
\begin{equation*}
S_{q}(z)=\frac{C_{q}}{x\left( z\right) }\sum_{s\geq 0}\frac{\left( q^{\beta
+1};q\right) _{s}\left( q^{1-N};q\right) _{s}(q^{-z};q)_{s}}{\left(
q^{1-N-\alpha };q\right) _{s}(q^{1-z};q)_{s}}\frac{q^{s}}{\left( q;q\right)
_{s}},\quad N+\alpha \notin\mathbb{N},
\end{equation*}
where 
\begin{equation*}
C_{q}=q^{\tilde{\upsilon}(\alpha ,\beta ,N)-1/2}\frac{\Gamma _{q}\left(
\beta +1\right) \Gamma _{q}\left( N+\alpha \right) }{\Gamma _{q}\left(
N\right) }.
\end{equation*}
Therefore 
\begin{equation}
S_{q}(z)=\frac{C_{q}}{x\left( z\right) }\,_{3}\varphi _{2}\left( 
\begin{array}{c|c}
q^{\beta +1},q^{1-N},q^{-z} &  \\ 
& q;q \\ 
q^{1-N-\alpha },q^{1-z} & 
\end{array}
\right) .  \label{Sq(z)H1}
\end{equation}
Finally, taking into account the relation 
\begin{equation*}
u_{0}^{q}=\frac{\left( q^{\alpha +\beta +2};q\right) _{N-1}}{\left(
q^{\alpha +1};q\right) _{N-1}}C_{q},
\end{equation*}
and using the transformation formula (see \cite[Theorem 12.4.2]{Ismail}) 
\begin{equation*}
\,_{3}\varphi _{2}\left( 
\begin{array}{c|c}
q^{-n},a,b &  \\ 
& q;q \\ 
c,d & 
\end{array}
\right) =\frac{b^{n}\left( d/b;q\right) _{n}}{\left( d;q\right) _{n}}
\,_{3}\varphi _{2}\left( 
\begin{array}{c|c}
q^{-n},b,c/a &  \\ 
& q;aq/d \\ 
c,q^{1-n}b/d & 
\end{array}
\right) ,
\end{equation*}
taking $n=N-1$, $a=q^{-z}$, $b=q^{\beta +1}$, $c=q^{1-z}$, and $d=q^{1-N-\alpha }$, relation (\ref{Sq(z)H1}) transforms into the following expression 
\begin{equation*}
S_{q}(z)=\frac{u_{0}^{q}}{x\left( z\right) }\,_{3}\varphi _{2}\left( 
\begin{array}{c|c}
q^{\beta +1},q^{1-N},q &  \\ 
& q;q^{N+\alpha -z} \\ 
q^{1-z},q^{\alpha +\beta +2} & 
\end{array}
\right),
\end{equation*}
where relation $\left( aq^{-n};q\right) _{n}=\left(
q/a;q\right) _{n}\left( -a\right) ^{n}q^{-n\left( n+1\right) /2}$ 
(see \cite[p. 304, (12.2.10)]{Ismail}) have been used.

\section{Conclusions and future directions}

Based on the accumulation of analytic and algebraic properties that
characterize the classical discrete orthogonal polynomials -often seemed to
be unrelated- the need for structure and classification of them has
constantly arisen as a central question in the orthogonal polynomial theory 
\cite{ga-ma-la,ma-la,me-al-ma,Niki}. In this paper we address this question
proving, in Theorem \ref{main-th}, that a sequence of $q$-polynomials
orthogonal with respect to $q$-moment functional (\ref{def-q-functional}) is
classical iff the associated Stieltjes function given in terms of $q$
-falling factorial basis satisfies a first order non-homogeneous $q$
-difference equation; the proof is given in a constructive way using a
theoretical background based on the theory of linear functional deeply
studied in \cite{maroni} and \cite{me-al-ma}. We show, in Theorem \ref
{main-th}, that the verification of the aforementioned difference equation
constitutes a new characterization of a discrete orthogonal polynomials on
the non-uniform lattice $x(s)=(q^{s}-1)/(q-1)$. However, more general
situations demand special attention. For instance, the $q$-orthogonal
polynomials on the lattice $x(s)=c_{1}q^{s}+c_{2}q^{-s}+c_{3}$, $(q\in 
\mathbb{R}^{+}\backslash {\{1\}})$, where the constants $c_{i}\in \mathbb{R}$
$(i=1,2,3)$ are constants independent of $s$, must be considered in the same
fashion as here. This paper outlines the important points and techniques to
be followed in such investigations aimed to characterize those families of $
q $-polynomials.

Finally, more general systems of $q$-moment functionals must be analyzed. In
this direction, the $q$-semiclassical orthogonal polynomials could be an
interesting challenge to be considered; in particular when Dirac masses are
added to $q$-moment functional. An analogous result to those given in
Theorem \ref{main-th} played a crucial role in the computation of the class
of the semiclassical functionals given by a perturbation via the addition of
Dirac masses (see \cite{garri-arvesu-paco}).

In closing, to the best of our knowledge, there is not in the literature any
explicit expression for the associated Stieltjes functions in the $q$
-falling factorial basis given in terms of hypergeometric functions.

\noindent\textbf{Acknowledgments} The research of the first author was
partially supported by the research grant MTM2009-12740-C03-01 of the
Ministerio de Educaci\'on y Ciencia of Spain and grant CC-G08-UC3M/ESP-4516 from Comunidad Aut\'o\-no\-ma de Madrid. We thank the reviewers for offering useful suggestions for improving the paper.

\end{document}